\documentclass[12pt]{article}
\usepackage{amssymb}
\usepackage{amsmath,amsfonts}

\def\ord{{\mbox{\ \rm ord\ }}}

\def\qq{{\mathbb Q}}
\def\cc{{\mathbb C}}
\def\zz{{\mathbb Z}}
\def\nn{{\mathbb N}}
\def\ff{{\mathbb F}}
\def\mod{{\mbox {\rm \ mod\ }}}
\def\ord{{\mbox {\rm ord}}}
\def\ind{{\mbox {\rm ind}_\omega}}

\def\ba{\begin{array}}
\def\ea{\end{array}}
\def\ds{\displaystyle}
\begin{document}

\begin{centering}

{\Large{\bf A cyclotomic approach to the solution of}

\Large{\bf  Waring's problem mod $p$}
}

\vspace*{ 0.5 cm}

M\'onica del Pilar Canales Ch.
\footnote{Supported by  DID--UACh, Grant  S--03--06}

\vspace*{ 0.4 cm}

Instituto de Matem\'aticas\\
Universidad Austral de Chile\\
Casilla 567,Valdivia

{\tt monicadelpilar@gmail.com}

\end{centering}

\vspace*{7 mm}


\textsc{Abstract}

Let $s_d(p,a) = \min \{k\ | \ a = \sum_{i=1}^{k}a_i^d, a_i\in \ff_p^*\}$
be the smallest number of $d$--th powers in the finite field $\ff_p$,
sufficient to represent the number $a\in \ff_p^*$.
Then 
$$g_d(p) = \max_{a\in\ff_p^*} s_d(p,a)$$ 
gives an  answer to Waring's Problem mod p.

We first introduce cyclotomic integers  $n(k,\nu)$, which then allow
to state and solve Waring's problem mod $p$ in terms of only the
cyclotomic numbers $(i,j)$ of order $d$.

We generalize the reciprocal of the Gaussian period equation $G(T)$ to a
$\cc$--differentiable function $I(T)\in\qq[[T]]$, which also satisfies
$I'(T)/I(T)\in\zz[[T]]$.
We show that and why $a\equiv -1\mod \ff_p^{*d}$ 
(the classical {\it Stufe}, if $d=2$)
behaves special: Here (and only here) $I(T)$ is in fact a polynomial
from $\zz[T]$, the reciprocal of the period polynomial.

We finish with  explicit calculations of $g_d(p)$ for the cases
$d=3$ and $d=4$, all primes $p$, using the known cyclotomic numbers
compiled by Dickson. 

\subsection*{1. Introduction}

Let $p > 2$ be a prime number, $d\geq 2$ a rational
integer, and let $s_d(p,a)$ be the least positive integer $s$ such that $a%
\in{\ff}_p^*$ is the sum of $s$ $d$--th powers in ${\ff}_p$, i.e.~
\[
s_d(p,a) = \min \{k\ |\ a=\sum_{i=1}^k a_i^d,\ a_i\in{\ff}_p^*\}.
\]

Since ${{\ff}_p^*}^d = {{\ff}_p^*}^{gcd(d,p-1)}$, it suffices to consider 
$d\ | \ p-1$. 

Let $f = (p-1)/d$ and let $\omega$ be a generator of ${{\ff}_p^*}$, fixed
from now on.
Clearly, it is enough to consider $s_d(p,a)$ only for the $d$ classes mod 
${{\ff}_p^*}^d$. 

Let $a\in {{\ff}_p^*}$ with $\alpha \equiv \ind (a) \mod d$ and 
let $\theta  \equiv \ind (-1)\mod d$, 
{\it i.e.} $\theta  =0$  if $f$ is even, and 
$\theta  = d/2$ if $f$ is odd.

In \cite{BC}\cite{C} we established, 
by considering the generating function
$$g(T) = \frac{1}{1-(\sum_{u\in \ff_p^{*d}}  \bar{X}^u)^T}
\in\left(K[X]/(X^p-1)\right)[[T]]$$
with $K = \qq(\zeta)$ the cyclotomic field given by $\zeta$ a
primitive $p$--th root of unity in $\cc$
that if
\[
N(k,a):=\# \{(u_{1},\dots ,u_{k})\in {{\ff}_{p}^{\ast }}^{d}\times \dots
\times {{\ff}_{p}^{\ast }}^{d}\ |\ a=u_{1}+\dots +u_{k}\}
\]
then
\[
N(k,a)=\frac{1}{p}\sum_{x=0}^{p-1}S(\zeta ^{x})^{k}\zeta ^{-ax}
= \frac{1}{p}\bigg[f^k + \sum_{\bar x\in{\ff}_p^*/{{\ff}_p^*}^d}
S(\zeta^{x})^k\cdot S(\zeta^{-ax})\bigg],
\]
where $S(\rho):=\sum_{u\in {{\ff}_{p}^{\ast }}^{d}}\rho ^{u}$.

Setting $i\equiv \ind(x)\mod d$
and $\alpha+\theta  \equiv \ind(-a)\mod d$,
we may write now 
\[
N(k,a)=\frac{1}{p}\bigg[f^{k}+\sum_{i=0}^{d-1}\eta _{i}^{k}\cdot \eta
_{i+\alpha +\theta  }\bigg]
\]
where
$$\eta _{i}:=S(\zeta ^{\omega ^{i}})=\zeta ^{\omega ^{i}}+\zeta ^{\omega
^{d+i}}+\zeta ^{\omega ^{2d+i}}+\dots +\zeta ^{\omega ^{(f-1)d+i}};0\leq
i\leq d-1$$
are the classical \textit{Gauss periods}, with minimal polynomial
over $\qq$ the so-called  {\it period polynomial of degree $d$}
\[
G(T)=\prod_{i=0}^{d-1}(T-\eta _{i})=\alpha _{d}+\alpha _{d-1}T+\dots +\alpha
_{2}T^{d-2}+T^{d-1}+T^{d}\in {\zz}[T],
\]
resolvent of the \textit{cyclotomic equation} $X^{p}-1=0$.

Hence, since  $s_d(p,a) = \min\{k \ | \ N(k,a)\neq 0\}$ we obtain:
$$s_d(p,a) = \min\{k \ | \ f^k+\sum_{i=0}^{d-1}
\eta_i^k\cdot \eta_{i+\alpha+\theta }\neq 0\}.$$

Our goal now is to determine $s_d(p,a)$, and thus $g_d(p)$, 
{\it only} using the {\it cyclotomic numbers of order $d$},
$$(i,j)\  :=\ \#\{(u,v);0\leq u,v\leq f-1\ |\ 1+\omega ^{du+i}\equiv \omega
^{dv+j}{\mbox {\rm \ mod\ }}p\},\ \ 0\leq i,j \leq d-1,$$
which have been extensively studied in the literature.

\subsection*{2. Cyclotomic Integers}

We call a polynomial expression  on the periods,  with integer
coefficients a {\it cyclotomic integer}
if it has an integer value.

Examples of this are the coefficients of
the period polynomial,
$$\alpha_k=s_k(\eta_0,\dots,\eta_{d-1}) =(-1)^k\sum_{0\leq
  i_1<\dots<i_k\leq d-1} 
\eta_{i_1}\ldots\eta_{i_k}\in\zz,$$
and its discriminant
$$D_d = \prod_{0\leq i<j\leq d-1} (\eta_i-\eta_j)^2\in \zz.$$

We study now, for all $k\in\nn$ and $0\leq \nu\leq d-1$, 
the  non trivial {\it cyclotomic integers}
$$n(k,\nu) = \sum_{i=0}^{d-1} \eta_i^k\cdot \eta_{i+\nu}\in\zz;\
0\leq \nu\leq d-1.$$

The theory of cyclotomy states the following formulae for the periods and
the cyclotomic numbers (see \cite{B}\cite{Di1}\cite{Di2}). 
For all $0\leq k,l\leq d-1$ it holds:

$
\begin{array}{ll}
(i) & \eta _{l}\eta _{l+k}=\sum_{h=0}^{d-1}(k,h)\eta _{l+h}+f\delta_{\theta k} \\
(ii) & \sum_{l=0}^{d-1}\eta _{l}\eta _{l+k}=p\delta_{\theta k}-f \\
(iii) & \sum_{h=0}^{d-1}(k,h)=f-\delta_{\theta k}
\end{array}
$

with Kronecker's $\delta_{ij}$.

The cyclotomic integers $n(k,\nu)$  assume values in
$\zz$ despite not being symmetric on the periods,
since they satisfy the following recurrence formula:

\textbf{Lemma 1.}\
{\it
Let $0\leq \nu \leq d-1$ and $k\geq 1$. Then
\[
n(k+1,\nu) = \sum_{l=0}^{d-1} (\nu,l) n(k,l) + f\delta_{\theta \nu} n(k-1,0)
\]
where $n(0,\nu)=-1$ and $n(1,\nu) = p\delta_{\theta\nu}-f$. 
}

\textit{Proof.}\ Clearly,
$$n(0,\nu)=\sum_{i=0}^{d-1}\eta_{i+\nu }=
\zeta+\zeta^2+\dots+\zeta^{p-1}=-1$$
and
$$n(1,\nu)=\sum_{i=0}^{d-1}\eta_{i+\nu }\eta_{i}=p\delta_{\theta \nu }-f$$

by $(ii)$.

Now by $(i)$, multiplying $\eta_l\eta_{l+\nu}$  by $\eta_l^k$ 
and adding over $l$ we get

$
\begin{array}{lllll}
n(k+1,\nu) & = & \sum_{l=0}^{d-1} \eta_{l+\nu}\eta_l^{k+1} &  &  \\
& = & \sum_{l=0}^{d-1}\left(\sum_{h=0}^{d-1}(\nu,h) \eta_{l+h}\eta_l^{k}+
f\delta_{\theta \nu}\eta_l^k\right) &  &  \\
& = & \sum_{h=0}^{d-1}(\nu,h)\sum_{l=0}^{d-1} \eta_{l+h}\eta_l^{k}+
f\delta_{\theta \nu} \sum_{l=0}^{d-1}\eta_l^k &  &  \\
& = & \sum_{h=0}^{d-1}(\nu,h)n(k,h) + f\delta_{\theta \nu} n(k-1,0). &
&
\hspace{25 mm}\Box
\end{array}
$

This turns out to be the key  to determine $N(k,a)$ in terms of
{\it only} the cyclotomic numbers in a remarkably simple  way, since we find


\textbf{Lemma 2.}\ \ 
{\it
Let $0 \leq \nu\leq d-1$. Then

$
\begin{array}{lll}
n(1,\nu)+f & = & p\delta_{\theta \nu } \\
n(2,\nu)+f^{2} & = & p(\nu ,\theta  ) \\
n(3,\nu)+f^{3} & = & p\sum_{i=0}^{d-1}(\nu ,i)(i,\theta 
)+f\delta_{\theta 0}[n(1,\nu)+f] \\
n(4,\nu)+f^{4} & = & p\sum_{i,j=0}^{d-1}(\nu ,i)(i,j)(j,\theta 
)+f\delta_{\theta 0}[n(2,\nu)+f^{2}]+f(0,\theta  )[n(1,\nu)+f].\\
\end{array}
$

and for $k \geq 5$

$
\begin{array}{lll}
n(k,\nu) + f^k & = & p\sum_{i_2,\dots,i_{k-1}=0}^{d-1}(\nu,i_2)%
\dots(i_{k-1},\theta )+ f\delta_{\theta 0}[n(k-2,\nu)+f^{k-2}] \\
& + & f(0,\theta )[n(k-3,\nu)+f^{k-3}] \\
& + & \sum_{j=4}^{k-1}\ f\ [\sum_{i_2,\dots,i_{j-2}=0}^{d-1}
(0,i_2)\dots(i_{j-2},\theta )][n(k-j,\nu)+f^{k-j}]
\end{array}
$ 
}

\textit{Proof.}\ 
The formulae for $k=1,2,3,4$ result by straightforward computation.
For higher $k$, we use induction on $k$. 
For all $0\leq l \leq d-1$, assume the formula for up to $k\geq 4$. This is
(the empty sum for $k=4$ is assumed null by convention)\\
$
\begin{array}{llll}
n(k,l) & = & p \sum_{i_2,\dots,i_{k-1}=0}^{d-1}(l,i_2)(i_2,i_3)
\dots(i_{k-1}, \theta )-f^k + f\delta_{\theta 0}[n(k-2,l)+f^{k-2}] &  \\
& + & f(0,\theta )[n(k-3,l)+f^{k-3}] &  \\
& + & \sum_{j=4}^{k-1}\ f\ [\sum_{i_2,\dots,i_{j-2}=0}^{d-1}
(0,i_2)\dots(i_{j-2},\theta )][n(k-j,l)+f^{k-j}]. &
\end{array}
$

Hence, by our previous lemma and the induction hypothesis, we have
\begin{eqnarray*}
n(k+1,\nu)
& = & \sum_{l=0}^{d-1} (\nu,l) n(k,l) + f\delta_{\theta \nu} n(k-1,0)\\
& = & \sum_{l=0}^{d-1} (\nu,l) \bigg\{ p\sum_{i_2,%
\dots,i_{k-1}=0}^{d-1}(l,i_2)(i_2,i_3)\dots(i_{k-1}, \theta )-f^k  \\
& + & f\delta_{\theta 0}[n(k-2,l)+f^{k-2}] +f(0,\theta )[n(k-3,l)+f^{k-3}]  \\
& + & \sum_{j=4}^{k-1}\ f\ [\sum_{i_2,\dots,i_{j-2}=0}^{d-1}
(0,i_2)\dots(i_{j-2},\theta )][n(k-j,l)+f^{k-j}]\bigg\}  \\
& + & f\delta_{\theta \nu} \bigg\{ p\sum_{i_2,\dots,i_{k-2}=0}^{d-1}(0,i_2)(i_2,i_3)
\dots(i_{k-2}, \theta )-f^{k-1}  \\
& + & f\delta_{\theta 0}[n(k-3,0)+f^{k-3}] +f(0,\theta )[n(k-4,0)+f^{k-4}]  \\
& + & \sum_{j=4}^{k-2}\ f\ [\sum_{i_2,\dots,i_{j-2}=0}^{d-1}
(0,i_2)\dots(i_{j-2},\theta )][n((k-1)-j,0)+f^{(k-1)-j}]\bigg\}. 
\end{eqnarray*}

Now, since
$p\delta_{\theta \nu }=n(1,\nu)+f$ and
$\sum_{l=0}^{d-1}(\nu,l)=f-\delta_{\theta \nu },$
and using that
\[
\sum_{l=0}^{d-1}(\nu ,l)n(k-j,l)=n(k+1)-j,\nu)-
f\delta_{\theta \nu}n((k-1)-j,0),
\]
the result follows as we may write
$\sum_{l=0}^{d-1}(\nu ,l)n(k-j,l)+(f-\delta_{\theta \nu })f^{k-j}$ as
$n((k+1)-j,\nu) +f^{(k+1)-j} 
 -f\delta_{\theta \nu }[n((k-1)-j,0)+f^{(k-1)-j}]$.\hfill $\Box $

\subsection*{3. The Cyclotomic Solution of Waring's Problem mod $p$}

We now may state the result that completes our study on higher levels 
and Waring's problem in $\ff_p$ in terms of cyclotomy.

\textbf{Theorem 1.}\ 
{\it
Let $p>2$ be a prime number and $d\geq 2$ an
integer with $p-1=df$. Let $\omega $ be a fixed generator of 
${\ff}_{p}^{\ast }$,
let $(i,j);0\leq i,j\leq d-1$ be the cyclotomic numbers of order $d$.
Also let $\theta  =0$ if $f$ even, and $\theta  =d/2$ if $f$ odd. 
Then, given $a\in {\ff}_{p}^{\ast }\backslash {\ff}_{p}^{{\ast }^d}$ with $\alpha \equiv \ind(a)
{\mbox {\rm \ mod\ }}d,$ we get
\[
s_{d}(p,a)=2\mbox{ \ if\ } (\alpha+\theta,\theta) \neq 0
\]
and otherwise
\[
s_{d}(p,a)=\min \{s\ |\ 
\exists\  0\leq i_2,\dots,i_{s-1}\leq d-1\colon
 (\alpha+\theta  ,i_{2})(i_{2},i_{3})\dots (i_{s-1},\theta )\neq 0 \}.
\]
}

\textit{Proof.}\ 
Since $a\in F_{p}^*\backslash F_{p}^{*^d}$ and 
$\alpha+\theta  \equiv \ind(-a)\mod d$, 
we have
\[
s_{d}(p,a)=\min \{k\geq 2\ |\ n(k,\alpha +\theta )+f^{k}\neq 0\},
\]
and hence $n(l,\alpha +\theta )+f^{l}=0$, for all $l<s_{d}(p,a)=s$. Then by
Lemma 2 
\[
n(s,\alpha +\theta )+f^{s}=p\sum_{i_{2},\dots ,i_{s-1}=0}^{d-1}(\alpha
+\theta  ,i_{2})(i_{2},i_{3})\dots (i_{s-1},\theta  ).
\]
Thus clearly
\[
s_{d}(p,a)=2,\mbox{\rm \ for\ }(\alpha+\theta  ,\theta  )\neq 0
\]
and otherwise $s_{d}(p,a)$ is the least integer $s$ with $3\leq s \leq d$, 
such that for some $0\leq i_2,\dots,i_{s-1}\leq d-1$, we have
a nonvanishing consecutive product of $s-1$  cyclotomic numbers
of the form
$ (\alpha +\theta ,i_{2})(i_{2},i_{3})\dots (i_{s-1},\theta  )\neq 0$.
\hfill $\Box $

Hence we obtain a solution of Waring's problem in $\ff_p$ 
via cyclotomy as:

\textbf{Theorem 2.}\ 
{\it
Let $p,d,f,\omega,\theta ,a,\alpha=\ind(a)$ 
and the cyclotomic numbers $(i,j)$ be as in Theorem $2$.
We define a matrix $M=(m_{ij})_{0\leq i,j\leq d-1}$ by\\
$m_{ij}= \left\{
\ba{cl}
0,&\mbox{\rm if\ } (i,j) = 0,\\
1,&\mbox{\rm otherwise,}\\
\ea
\right.$
and we denote  its $n$--th power as $(m_{ij}^{(n)}) := M^n$.

Then
\[
g_d(p) = \max_{0\leq \alpha\leq d-1}\min\{ s \ |\ 
m_{(\alpha+\theta )\theta }^{(s-1)} \neq 0\}.
\]
}

\textit{Proof.}\ 
By Theorem 1, $s_{d}(p,a)$ is the least integer $s$ with $2\leq s \leq d$ 
such that $(\alpha+\theta,\theta)\neq 0$
or for some $0\leq i_2,\dots,i_{s-1}\leq d-1$ we have
$ (\alpha +\theta ,i_{2})(i_{2},i_{3})\dots (i_{s-1},\theta  )\neq 0$.
Also, $\ds g_d(p) = \max_{0\leq \alpha\leq d-1}\{s_d(p,a)\}$.

Thus, since 
$$m^{(n)}_{(\alpha+\theta)\theta } = \sum^{d-1}_{i_2,i_3,\dots,i_{n}=0}
m_{(\alpha+\theta )i_2}\cdot m_{i_2i_3}\cdot\ \dots\  \cdot m_{i_n\theta }$$
is the entry at $(\alpha+\theta ,\theta )$ of the $n$--th power of the matrix 
$M$, where $m_{ij}\neq 0$ iff $(i,j)\neq 0$, the result follows.
\hfill$\Box$


\subsection*{4. On the Generalization of Theorem 1 in \cite{BC}}

In \cite{BC}, we only consider the case  $a\equiv -1 \mod \ff_p^{*d}$,
{\it i.e.} $\alpha+\theta  \equiv 0\mod d$, finding 
$s_d(p,-1)$ in terms of the coefficients 
$\alpha_k; 2\leq k \leq d$ of the period polynomial.

Now if $\alpha +\theta  \not\equiv 0\mod d$, we can
generalize Theorem 1 in \cite{BC} as follows:

{\bf Theorem 3.}\ 
{\it 
Let $p>2$  be a prime number and $d\geq 2$ an integer with $d\ | \ p-1$.
Let $\omega$ be a fixed generator of $\ff_p^*$ and
let $\eta_i; 0\leq i\leq d-1$ be the Gaussian periods. 
Then if $a\in\ff_p^*\backslash \ff_p^{*d}$ with 
$\ind(-a) \equiv \alpha + \theta  \mod d$, 
we have
$$s_d(p,a) = \ord_T\left(\frac{1}{1-fT}
-\frac{I_{\alpha+\theta }(T)'}{I_{\alpha+\theta }(T)}\right)$$
where 
$$I_{\alpha+\theta }(T)=\prod_{i=0}^{d-1}
(1-\eta_iT)^{\left(\frac{\eta_{i+\alpha+\theta }}{\eta_i}\right)}\in\qq[[T]]$$
is a complex differentiable function of $T$ and $\ord_T$ is the 
usual valuation in $\zz[[T]]$.
}

{\it Proof:}\ 
We have
$$s_d(p,a)=\min\{k\ | \ N(k,a) \neq 0\} = \ord_T(\sum_{k=0}^\infty
N(k,a)T^k)$$
where 
$$N(k,a) 
= \frac{1}{p}[f^k + n(k,\alpha+\theta)]
= \frac{1}{p}[f^k + \sum_{i=0}^{d-1} \eta_i^k\cdot \eta_{i+\alpha+\theta}].$$
Thus formally
\begin{eqnarray*}
\sum_{k=0}^\infty N(k,a)T^k
 &=& \frac{1}{p}[\sum_{k=0}^\infty f^kT^k +  
\sum_{k=0}^\infty  n(k,\alpha+\theta)T^k]\\
 &=& \frac{1}{p}[\sum_{k=0}^\infty f^kT^k + \sum_{k=0}^\infty
(\sum_{i=0}^{d-1}   \eta_i^k\cdot \eta_{i+\alpha+\theta})T^k]\\
 &=& \frac{1}{p}[\frac{1}{1- fT} + \sum_{i=0}^{d-1} 
(\sum_{k=0}^\infty  \eta_i^k\cdot \eta_{i+\alpha+\theta}T^k)]\\
 &=& \frac{1}{p}[\frac{1}{1- fT} + \sum_{i=0}^{d-1} 
\frac{ \eta_{i+\alpha+\theta}}{1-  \eta_i T}].\\
\end{eqnarray*}

Now, considering 
$\ds\frac{\eta_{i+j}}{1-\eta_iT}$ as a  complex differentiable 
function of $T$  for  all $0\leq i,j \leq d-1$, and recalling that 
$\ds\frac{\eta_{i+j}}{\eta_i}\in \cc\backslash\{0\}$, we have\\
$
\ba{rcl}
\ds\sum_{i=0}^{d-1} \frac{\eta_{i+j}}{1-\eta_iT}
&=&\ds-\sum_{i=0}^{d-1} \frac{\eta_{i+j}}{\eta_i}\cdot
\frac{(1-\eta_{i}T)'}{(1-\eta_iT)}
\\
&=&\ds-\sum_{i=0}^{d-1} \frac{\eta_{i+j}}{\eta_i}[\log(1-\eta_{i}T)]'\\
&=&\ds-[\sum_{i=0}^{d-1} \log((1-\eta_{i}T)^{\frac{\eta_{i+j}}{\eta_i}})]'\\
&=&\ds-[\log(\prod_{i=0}^{d-1} (1-\eta_{i}T)^{\frac{\eta_{i+j}}{\eta_i}})]'\\
&=&\ds-[\log I_j(T)]'\\
&=&\ds-\frac{I_j(T)'}{I_j(T)} \ \ \in \zz[[T]],
\ea$

where $I_j(T) = \prod_{i=0}^{d-1} (1-\eta_iT)^{\frac{\eta_{i+j}}{\eta_i}}$
has a  power series expansion $I_j(T) = \sum_{k=0}^\infty c_kT^k\in\cc[[T]]$, 
where $c_k = - \frac{1}{k}\sum_{l=0}^{k-1}c_l\cdot n(k-1-l,j)$ and 
$c_0=1$. Thus,  we recursively obtain 
$k!c_k\in\zz,\forall k\geq 0$, and  hence
$I_j(T),I_j(T)'\in\qq[[T]]$.~\hfill~$\Box$ 


{\it Remark:}\ This  finally shows the class of $-1$ 
to be special since 
$$I_{\alpha+\theta }(T)\in\qq[T] \Leftrightarrow 
\frac{\eta_{i+\alpha+\theta }}{\eta_i}\in\nn,
\ \forall\ 0\leq i\leq d-1.$$

This is, iff  $\alpha+\theta =0$ and hence $a\equiv -1 \mod \ff_p^{*d}$.
In this case, 
$$I_0(T) = \prod_{i=0}^{d-1} (1-\eta_{i}T)= T^dG(T^{-1})$$ 
is the reciprocal of the Gauss period polynomial 
and we recover Theorem 1 of \cite{BC}.

\subsection*{5. Explicit Numerical Results}

We state the complete results for $d=3$ and $d=4$, for all primes $p$,
following  \cite{Di1}\cite{Di2}:

\textbf{Theorem 4.}\ \textit{Let $p=3f+1$ be a prime number with $%
4p=L^{2}+27M^{2}$ and 
$L\equiv 1{\mbox {\rm \ mod\ }}3$. Then
\[
g_3(p)=\bigg\{
\begin{array}{ll}
3, & \mbox{ \ if \ }p=7, \\
2, & \mbox{\ otherwise.}
\end{array}
\]
}

\textit{Proof.}\ Since $f$ is even, we have $(h,k)=(k,h)$ , and it is known
by \cite{Di1} that\\
\[
18(0,1)=2p-4-L+9M\text{ and }\\
18(0,2)=2p-4-L-9M.
\]

Then by Theorem 1, with $\theta  =0$ and the sign of $M$ depending on the
choice of the generator $\omega$, we have

\[
s_{3}(p,\omega)=\bigg\{
\begin{array}{ll}
2, & \mbox{\rm \ if \ }\left( 1,0\right)\neq 0
\mbox{\rm \ \ i.e. \ }2p\neq 4  + L- 9M \\
3, & \mbox{\rm \ otherwise \ }
\end{array}
\]

and

\textit{
\[
s_{3}(p,\omega^{2})=\bigg\{
\begin{array}{ll}
2, & \mbox{\rm \ if \ }\left( 2,0\right) \neq 0\mbox{\rm \ \ i.e. \ }2p\neq 4%
+L+9M \\
3, & \mbox{\rm \ otherwise. \ }
\end{array}
\]
}
\newpage
Thus

\textit{
\[
g_3(p)=\bigg\{
\begin{array}{ll}
2, & \mbox{\rm \ if \ }2p\neq 4+L\pm 9M, \\
3, & \mbox{\rm  \ otherwise. \ }
\end{array}
\]
}

Now, $g_3(p) = 3
\Longleftrightarrow \exists \alpha\in\{1,2\}$ with $(\alpha,0)=0
\Longleftrightarrow 
4p=L^{2}+27M^{2}=8+2L\pm 18M
\Longleftrightarrow L=M=1
\Longleftrightarrow p=7$.
\hfill $\Box $

{\bf Theorem 5.}
\ \textit{Let $p=4f+1$ be a prime number with 
$p=x^{2}+4y^{2}$ and 
$x\equiv 1{\mbox {\rm \ mod\ }}4$. Then
\[
g_4(p)=\Bigg\{
\begin{array}{ll}
4, & \mbox{ \ if \ }p=5, \\
3, & \mbox{ \ if \ }p=13, 17, 29, \\
2, & \mbox{ \ otherwise.\ }
\end{array}
\]
}

\textit{Proof.} \ By Theorem 1, with  
$\theta  =\bigg\{$
\begin{tabular}{cc}
0, & if $f$ even,\\
$d/2$, & if $f$ odd,
\end{tabular}
and the sign of $y$\\ 
depending on the choice of the generator $\omega$,  we have

$s_4(p,\omega^\alpha) =\left\{
\begin{array}{lll}
1,& \mbox{ \rm  if \ } \alpha=0 \\
2,& \mbox{ \rm  if \ } \alpha\neq 0, (\alpha+\theta ,\theta )\neq 0 \\
3,& \mbox{ \rm  if \ } \alpha\neq 0, (\alpha+\theta ,\theta )= 0\\
&\mbox{\ \rm  and \ } (\alpha+\theta ,i)(i,\theta )\neq 0 
\mbox{ \rm  for some \ } 0\leq i\leq 3\\
4,&\mbox{ \rm  otherwise \ }
\end{array}
\right.
$

where  by \cite{Di1} we may find the cyclotomic numbers in terms of the 
representation of $p$.

If $f$ is even:

$16(0,0) = p-11-6x$\\
$16(0,1) = p-3+2x+8y$\\
$16(0,2) = p-3+2x$\\
$16(0,3) = p-3+2x-8y$\\
$16(1,2) = p+1-2x$

and

$(1,1) = (0,3)$, 
$(1,3) = (2,3) = (1,2)$, 
$(2,2) = (0,2)$, 
$(3,3) = (0,1)$, with 
$(i,j) = (j,i)$. 

\newpage

If $f$ is odd:

\nopagebreak

$16(0,0) = p-7+2x$\\
$16(0,1) = p+1+2x-8y$\\
$16(0,2) = p+1-6x$\\
$16(0,3) = p+1+2x+8y$\\
$16(1,0) = p-3-2x$

and

$(1,1) = (2,1) = (2,3) = (3,0) = (3,3) = (1,0)$, 
$(1,2) = (3,1) = (0,3)$, 
$(1,3) = (3,2) = (0,1)$, 
$(2,0) = (2,2) = (0,0)$. 

Thus, we find $g_4(p) > 2$ if $p=x^2+4y^2,$ with $x\equiv 1 \mod 4$,
satisfies one of the following  diophantine equations:

$f$ even:

$\ba{lll}
(\alpha=1)& x^2+4y^2+2x+8y &= 3\\
(\alpha=2)& x^2+4y^2+2x  &= 3\\
(\alpha=3)& x^2+4y^2+2x-8y& = 3
\ea$

$f$ odd:

$\ba{lll}
(\alpha=1)  &   x^2+4y^2+2x-8y& = -1 \\
(\alpha=2)  &   x^2+4y^2-6x   & = -1 \\
(\alpha=3)  &   x^2+4y^2+2x+8y& = -1
\ea$ 

All these are equations of the form $(x+a)^2 + 4(y+b)^2=c$ and one easily 
finds that the only solutions give  $p= 5,13,17,$ and 29.

Now, checking for these primes the equations for $g_4(p)=4$ we find only
$g_4(5)=4$,  and thus $g_4(13) = g_4(17) = g_4(29) = 3$.\hfill $\Box$

We may obtain complete solutions for Waring's problem mod $p$
and thus 
Waring's problem mod $n$ (see \cite{SAMM1}\cite{SAMM2}), for all $d\geq 3$ 
for which the cyclotomic numbers are known or may be found in terms of 
the representations of multiples of $p$ by binary quadratic forms.
Clearly, for $d> 4$ much effort is needed to obtain the $d^2$ cyclotomic 
constants, and other representations of multiples of $p$ by quadratic forms 
and the study of different cases is necessary
(see~\cite{BC}\cite{C}\cite{Di1}\cite{Di2}),
but since $s_d(p,-1)$ and $g_d(p)$ are still open problems for 
$3d+1\leq p < (d-1)^4$, our work  seems  to give the only 
complete theoretical result on the subject.

\subsection*{Conclusion}

We introduced the concept of cyclotomic integers and gave some
non--tri\-vial examples, the $n(k,\nu)$, which allowed us to solve 
the modular Waring's  Problem using only the classical
cyclotomic numbers.

We saw that the analytic function 
$I_{\alpha+\theta}(T)$ is a formal power series with coefficients in
$\qq$,
and that for $a\equiv -1 \mod \ff_p^{*d}$,
{\it i.e.} $\alpha = \theta $, 
in fact $I_{\alpha+\theta}(T)=I_0(T)$ 
is a polynomial in $\zz[T]$, the reciprocal of the
Gauss period equation.

We finished with two examples of  explicit calculations of $g_d(p)$ 
for $d=3$ and $d=4$, all primes $p$.

\end{document}